\documentclass{amsart}
\usepackage[T1]{fontenc}
\usepackage[latin1]{inputenc}
\makeatletter

\theoremstyle{plain}    
\newtheorem{thm}{Theorem}[section]
\numberwithin{equation}{section} 
\numberwithin{figure}{section} 
\theoremstyle{plain}    
\newtheorem{cor}[thm]{Corollary} 
\theoremstyle{plain}    
\newtheorem{lem}[thm]{Lemma} 
\theoremstyle{plain}    
\newtheorem{prop}[thm]{Proposition} 
\theoremstyle{definition}
\newtheorem{defn}[thm]{Definition}

\usepackage{hyperref}
\usepackage[arrow,matrix,curve,ps]{xy}
\CompileMatrices

\usepackage{epsfig}

\renewcommand{\O}{\mbox{$\mathcal{O}$}}
\renewcommand{\P}{\mbox{$\mathbb{P}$}}
\newcommand{\Q}{\mbox{$\mathbb{Q}$}}
\newcommand{\C}{\mbox{$\mathbb{C}$}}
\newcommand{\Z}{\mbox{$\mathbb{Z}$}}

\DeclareMathOperator{\Chow}{Chow}

\DeclareMathOperator{\Div}{Div}

\DeclareMathOperator{\Image}{Image}

\DeclareMathOperator{\mult}{mult}
\DeclareMathOperator{\Pic}{Pic}

\DeclareMathOperator{\Reg}{Reg}
\DeclareMathOperator{\Sing}{Sing}

\DeclareMathOperator{\Spec}{Spec}

\makeatother

\begin{document}

\title{Projective bundles of singular plane cubics }

\author{Stefan Kebekus}

\date{\today}

\address{Institut für Mathematik, Universität Bayreuth, 95440
  Bayreuth, Germany}

\email{stefan.kebekus@uni-bayreuth.de}

\urladdr{http://btm8x5.mat.uni-bayreuth.de/\protect\( \sim \protect \)kebekus}

\thanks{The author gratefully acknowledges support by a
  Forschungsstipendium of the Deutsche Forschungsgemeinschaft.}

\begin{abstract}
  Classification theory and the study of projective varieties which
  are covered by rational curves of minimal degrees naturally leads to
  the study of families of singular rational curves. Since families of
  arbitrarily singular curves are hard to handle, it has been shown in
  \cite{Keb00a} that there exists a partial resolution of
  singularities which transforms a bundle of possibly badly singular
  curves into a bundle of nodal and cuspidal plane cubics.
  
  In cases which are of interest for classification theory, the total
  spaces of these bundles will clearly be projective. It is, however,
  generally false that an arbitrary bundle of plane cubics is globally
  projective. For that reason the question of projectivity seems to be
  of interest, and the present work gives a characterization of the
  projective bundles.
\end{abstract}
\maketitle
\tableofcontents

\section{Introduction}

Let $C$ be a smooth algebraic curve and $\pi :X\to C$ be a morphism
from a (singular) algebraic variety to $C$. We assume that every fiber
of $\pi$ is isomorphic to an irreducible and reduced singular plane
cubic. Omitting the word ``integral'' for brevity, we call this setup
a ``bundle of singular plane cubics''. Although $C$ can always be
covered by open subsets $U_{\alpha}$ such that $X_{\alpha} := \pi
^{-1}(U_{\alpha})$ can be identified with a family of cubic curves in
$\P_2\times U_{\alpha}$, it is generally not true that $X$ can be
embedded into a $\P_2$-bundle over $C$. In fact, $X$ even need not be
projective.  The aim of this paper is to characterize those bundles
which \emph{are} projective.

This problem arises naturally in the study of projective varieties
which are covered by a family of rational curves of minimal degrees:
if we are given a projective variety $V$ and a proper subvariety
$$
H\subset \Chow (V)
$$
parameterizing a covering family of rational curves of minimal
degree, the subfamily $H^{\Sing} \subset H$ parameterizing singular
rational curves is of greatest interest. It is conjectured that $\dim
H^{\Sing }<\dim V-1$. In our previous paper \cite{Keb00a} we gave
bounds on the dimension of $H^{\Sing}$ using the following line of
argumentation. First, we chose a general point $x\in V$, considered
the subfamily $H_x^{\Sing} \subset H^{\Sing}$ of singular curves which
contain $x$ and constructed a diagram as follows:
$$
\xymatrix{ {\tilde U} \ar[rr]^{\txt{\scriptsize normalization}}
  \ar@/_/[rrd]_{\txt{\scriptsize $\P_1$-bundle}} & & {U'}
  \ar[rrr]^{\txt{\scriptsize finite morphism}}
  \ar[d]^{\txt{\scriptsize bundle of singular\\ \scriptsize plane
      cubics}} & & & {U} \ar[d]^{\pi} \\ & & {\tilde H}
  \ar[rrr]_{\txt{\scriptsize finite cover and \\ \scriptsize
      normalization}} & & & {H^{\Sing}_x}} 
$$
Where $U\subset \Chow (X)\times X$ is the universal family and
fibers of $\pi $ are irreducible and generically reduced singular
rational curves.

Secondly, we noted that the family $U'$ comes from the universal
family over the Chow-variety and is therefore projective. We were able
to show that this is possible only if the parameter space
$H_x^{\Sing}$ is either finite or if it is 1-dimensional and
parameterizes both nodal and cuspidal curves.  The argumentation
involved an analysis of the intersection of certain divisors in
$\tilde{U}$. A complete description of projective bundles, which is
somewhat more delicate and not numerical in nature, has not been given
in \cite{Keb00a}. In order to complete the picture we discuss it
here. It is hoped that these results will be useful in the further
study of rational curves on projective varieties.

Throughout the present paper we work over the field $\C$ of complex
numbers and use the standard language of algebraic geometry as
introduced in \cite{Ha77}.

\subsubsection*{Acknowledgement}

The results of this paper were worked out while the author enjoyed the
hospitality of RIMS in Kyoto and KIAS in Seoul. The author is grateful
to Y.~Miyaoka for the invitation to RIMS and to J.-M.~Hwang for the
invitation to KIAS. He would also like to thank S.~Helmke for a number
of helpful discussions.

\section{Ruled surfaces and elementary transformations}

The main results of this paper characterize projective bundles of
singular plane cubics by describing their normalizations, which are
ruled surfaces. Thus, before stating the results in
section~\ref{sec:main_results}, it seems advisable to recall some
elementary facts about the normalization morphism and about elementary
transformations between ruled surfaces.

\subsection{Reduction to ruled surfaces}

\label{sec:Reduct_to_ruled}
As a first step in the reduction of the characterization problem, we
note that to give a bundle of singular plane cubics over a smooth
curve $C$, it is equivalent to give a ruled surface $\tilde{X}$ over
$C$ and a double section $\tilde{\sigma} \subset \tilde{X}$.  Indeed,
if a bundle $X$ of singular plane cubics is given, we know from
\cite[thm.~II.2.8]{K96} that its normalization will be a
$\P_1$-bundle. The scheme-theoretic preimage of the (reduced) singular
locus will be a double section. On the other hand, if $\tilde{\pi} :
\tilde{X}\to C$ is a $\P_1$-bundle containing a double section
$\tilde{\sigma}$, we can construct a diagram 
$$
\xymatrix{ {\tilde X} \ar[rr]^{\gamma}_{\txt{\scriptsize
      identification}} \ar[rd]_{\txt{\scriptsize $\tilde \pi$ \\
      \scriptsize
      $\P_1$-bundle}} & & {X} \ar[ld]^{\txt{ \scriptsize $\pi$ \\
      \scriptsize bundle of singular\\ \scriptsize plane cubics}} \\ &
  C} 
$$
as follows. Find a cover by open subsets $U_{\alpha }\subset C$ so
that we can identify $\tilde{\pi }^{-1}(U_{\alpha })\cong \P_1 \times
U_{\alpha }$ in a way that enables us to write
$$
\tilde{\sigma }\cap \tilde{\pi }^{-1}(U_{\alpha })=\left\{
  ([y_{0}:y_{1}],x)\in \P _{1}\times U_{\alpha }\, |\,
  y_{0}^{2}=g(x)y^{2}_{1}\right\}
$$
where $g\in \O (U_{\alpha })$. The identification morphism $\gamma
_{\alpha }=\gamma |_{\tilde{\pi }^{-1}(U_{\alpha })}$ is then locally
given as
$$
\begin{array}{cccc}
  \gamma _{\alpha }: & \P _{1}\times U_{\alpha } & \to  & \P _{2}\times
  U_{\alpha }\\ 
  & ([y_{0}:y_{1}],x) & \mapsto  & \left(
    [y^{2}_{0}y_{1}-g(x)y^{3}_{1}:
    y^{3}_{0}-g(x)y_{0}y^{2}_{1}:y^{3}_{1})],x\right)
\end{array}
$$
An elementary calculation shows that the image of $\gamma_{\alpha}$
has the structure of a bundle of singular plane cubics and that the
local morphisms $\gamma_{\alpha}$ glue together to give a global one.
More precisely, for a point $\mu \in C$ with fiber $X_{\mu}:=\pi
^{-1}(\mu)$, the fiber is nodal if $\tilde{\sigma}$ intersects
$\tilde{X}_{\mu}:=\tilde{\pi}^{-1}(\mu)$ in two distinct points and
cuspidal if it intersects in a double point. Note that the gluing
morphisms do not in general come from automorphisms of $\P_2$. For
that reason $\tilde{X}$ can in general \emph{not} be embedded into a
$\P_2$-bundle over $C$.

\subsection{Elementary transformations}

The primary tool in the discussion of ruled surfaces will be the
``elementary transformation'' which is a birational map between ruled
surfaces. We refer to \cite[V.5.7.1]{Ha77} for the definition and a
brief discussion of these maps and the associated terminology.

If $\pi :Y\to C$ is a ruled surface and $(\sigma_i,D_i)_{i=1\ldots n}$
is a collection of sections $\sigma_i\subset Y$ and effective divisors
$D_i\in \Div (C)$ such that the supports $|D_i|$ are mutually
disjoint, we can inductively define a birational map between ruled
surfaces
$$
elt_{(\sigma _{i},D_{i})_{i=1\ldots n}}:Y\dasharrow \tilde{Y}
$$
as follows. Choose an index $j$, choose a closed point $\mu \in
|D_{j}|$ and perform an elementary transformation with center $\pi
^{-1}(\mu )\cap \sigma_j$. Replace the $\sigma_i$ with their strict
transforms, replace $D_i$ with $D_i-\delta_{ij}\mu $, where $\delta $
is the Kronecker symbol, and start anew until all $D_i$ are zero. It
follows directly from the construction of the elementary
transformation that the target variety $\tilde{Y}$ as well as the
resulting birational map are independent of the choices made.

The inverse of an elementary transformation can again be written as an
elementary transformation and the following lemma shows a way to write
down the inverse transformation locally. The proof is very elementary
and therefore omitted here.

\begin{lem}\label{lem:inverse_trafo}
  Let $\pi :Y\to C$ be a ruled surface ($C$ not necessarily compact)
  and let $(\sigma_i)_{i=1\ldots 3}$ be sections. Let $D_{1}\in \Div
  (C)$ be an effective divisor and consider the birational map
  $elt_{(\sigma _{1},D_{1})}:Y\dasharrow \tilde{Y}$.  If
  $\tilde{\sigma }_{i}\subset \tilde{Y}$ are the strict transforms of
  the $\sigma _{i}$, then the following holds:
  \begin{enumerate}
  \item If $\sigma_1$, $\sigma_2$ and $\sigma_3$ are mutually
    disjoint, then $D_1$ can be written as
    $$
    D_1=\sum _{p\in \tilde{\sigma }_2\cap \tilde{\sigma }_3}
    \mult_p(\tilde{\sigma }_2,\tilde{\sigma }_3)\cdot \tilde{\pi }(p)
    $$
    Here $\mult _{p}(\tilde{\sigma }_{2},\tilde{\sigma }_{3})$
    denotes the local intersection number of $\tilde{\sigma }_{2}$ and
    $\tilde{\sigma }_{3}$ at $p$.
  
  \item If $\sigma_1$ and $\sigma_2$ are disjoint, then the inverse
    map $elt_{(\sigma_1,D_1)}^{-1} : \tilde{Y}\dasharrow Y$ is given
    as $elt_{(\sigma_1,D_1)}^{-1}=elt_{(\tilde{\sigma }_2,D_1)}$.
  \end{enumerate}
  \hfill\qed
\end{lem}

A repeated application of lemma~\ref{lem:inverse_trafo}.(2) allows us
to write the inverse transformation in a more complicated situation.
Again we omit the proof.

\begin{cor}\label{cor:inverse_trafo}
  Let $\pi :Y\to C$ be a ruled surface and let $(\sigma_i)_{i=1\ldots
    n}$ be sections and $D_i \in \Div(C)$ be effective divisors with
  disjoint supports. Assume that for every index $i$ and every point
  $\mu \in |D_i|$, there exists a unique index $j$ such that $\pi
  (D_i\cap D_j)\not \ni \mu $. Consider the birational map
  $elt_{(\sigma_i,D_i)_{i=1\ldots n}}:Y\dasharrow \tilde{Y}$. If
  $\tilde{\sigma}_i \subset \tilde{Y}$ are the strict transforms of
  the $\sigma_i$, and if we set
  $$
  \tilde{D}_i := \sum _{j\not =i,\, \mu \in |D_j|\setminus
    \pi(\sigma_i \cap \sigma_j)}\mult_{\mu} (D_j)\cdot \mu ,
  $$
  then the inverse map is given as $elt_{(\sigma_i, D_i)_{i=1\ldots
      n}}^{-1}=elt_{(\tilde{\sigma}_i,\tilde{D}_i)_{i=1\ldots n}}$.
  \hfill\qed
\end{cor}

\section{Characterization of projective bundles}
\label{sec:main_results}

The following two theorems are the main results of this paper.
Theorem~\ref{thm:construction} gives a construction which yields
examples for projective and non-projective bundles of singular plane
cubics. Theorem~\ref{thm:main_characterization} shows that ---after
finite base change, if necessary--- all projective bundles of singular
plane cubics can be constructed by this method.

\begin{thm}\label{thm:construction}
  Let $C$ be a smooth curve, let $n$ a positive integer,
  $(D_i)_{i=1\ldots n}\in \Div (C)$ arbitrary disjoint effective
  divisors and $\sigma_0$, $(\sigma_i)_{i=1\ldots n}$ and
  $\sigma_{\infty }$ arbitrary distinct fibers of the projection
  $\P_1\times C\to \P_1$.  Construct a bundle $X$ of singular plane
  cubics as follows. 
  \begin{equation}\label{eq:backtransform}
    \xymatrix{ {\P_1 \times C} \ar@{-->}[rrr]^{elt_{(\sigma_i,
          D_i)_{i=1\ldots n}}}_{\txt{\scriptsize elementary\\
          \scriptsize transformations}} \ar[d]_{\txt{\scriptsize
          projection}}^{\pi_2} & & & {\tilde X}
      \ar[rr]^{\gamma_{(\tilde \sigma_0, \tilde
          \sigma_\infty)}}_{\txt{\scriptsize identification}}
      \ar[d]_{\tilde \pi} & & X \ar[d]^{\txt{\scriptsize bundle of
          singular\\ \scriptsize plane cubics}}_{\pi} \\ C
      \ar@{=}[rrr] & & & C \ar@{=}[rr]& & C} 
  \end{equation} 
  where $\gamma_{(\tilde{\sigma}_0, \tilde{\sigma}_{\infty})}$ is the
  identification morphism described in
  section~\ref{sec:Reduct_to_ruled}. The bundle $X$ is projective if
  and only if there exists a coordinate on $\P _{1}$ such that
  $\sigma_0 =\{[0:1]\}\times C$, $\sigma _{\infty }=\{[1:0]\}\times C$
  and $\sigma_i=\{[\xi_i:1]\}\times C$ where the $\xi_i$ are roots of
  unity.
\end{thm}
  
\begin{thm}\label{thm:main_characterization}
  Let $\pi :X'\to C'$ be a projective bundle of singular plane cubics
  over a smooth curve $C'$ and assume that nodal curves occur as
  fibers. Then after finite base change, if necessary, $X'$ can be
  constructed by the method described in
  theorem~\ref{thm:construction} above. More precisely, there exists a
  finite morphism $\tau :C\to C'$ between smooth curves, there exist
  effective divisors $D_i\in \Div (C)$ and disjoint sections
  $\sigma_0$, $\sigma_\infty$ and $(\sigma_i)_{i=1\ldots n}\subset
  \P_1\times C$ such that $X'$ can be constructed in the following
  manner. 
  \begin{equation}
    \xymatrix{ {\P_1 \times C}
      \ar@{-->}[rrr]^{elt_{(\sigma_i, D_i)_{i=1\ldots
            n}}}_{\txt{\scriptsize elementary\\ \scriptsize
          transformations}} \ar[d]_{\txt{\scriptsize
          projection}}^{\pi_2} & & & {\tilde X}
      \ar[rr]^{\gamma_{(\tilde \sigma_0, \tilde
          \sigma_\infty)}}_{\txt{\scriptsize identification}}
      \ar[d]_{\tilde \pi} & & X \ar[d]_{\pi} \ar[rr] && X'
      \ar[d]^{\txt{\scriptsize bundle of singular\\ \scriptsize plane
          cubics}}_{\pi'} \\ C \ar@{=}[rrr] & & & C \ar@{=}[rr]& &
      C\ar[rr]_{\txt{\scriptsize finite base change}}^{\tau} && C'}
  \end{equation} 
  Here $X$ denotes the fibered product $X:=X'\times_{C'}C$.
\end{thm}

\section{Projective Bundles}

The present section is concerned with an investigation of the special
geometry of projective bundles of singular plane cubics. The results
will later be used in sections~\ref{sec:Proof1} and \ref{sec:Proof2}
to prove the main theorems.  Throughout this section we assume that
$\pi :X\to C$ is a projective bundle of singular plane cubics over a
smooth curve $C$ and that $L\in \Pic (X)$ is an ample line bundle.
Furthermore we assume that both nodal and cuspidal cubics occur as
fibers of $\pi $. In this setup deformation theory tells us that the
generic fiber has a nodal singularity and that there are only finitely
many fibers with cusps.

Let $\eta :\tilde{X}\to X$ be the normalization. By
\cite[thm.~II.2.8]{K96}, the variety $\tilde{X}$ is a $\P_1$-bundle
over $C$.

\subsection{\protect$L\protect$-osculating points}

The restriction of the line bundle $L$ to a fiber with nodal
singularity defines a number of points which we call
``$L$-osculating''. More precisely, we use the following definition.

\begin{defn}\label{not:special_points}
  Let $C^N$ be a nodal plane cubic, and $H\in \Pic ^{k}(C^N)$ be a
  line bundle of degree $k>0$. We call a smooth point $\sigma \in
  C_{\Reg }^N$ an ``$H$-osculating point'' if $\O (k\sigma )\cong H$.
\end{defn}

The following lemma shows how to calculate the $H$-osculating points
on a given curve in a particularly simple situation.

\begin{lem}\label{lem:how_to_calculate_special_pts}
  Let $C^{N}$ be a nodal plane cubic and $p\in C_{\Reg }^{N}$ be a
  smooth point. Fix an identification $\iota :\C ^{*}\to C_{\Reg}^N$
  such that $\iota^{-1}(p)=1$ and set
  $$
  (\sigma_i)_{i=1\ldots k}=\{\iota (\xi)\, |\, \xi^k=1\}.
  $$
  Then $\sigma_i$ are the osculating points for the line bundle
  $\O(kp)$. In particular, there exist exactly $k$ osculating points
  for $\O(kp)$.
\end{lem}

\begin{proof}
  Recall from \cite[Ex. II.6.7]{Ha77} that the map $\iota$ defines a
  group morphism
  $$
  \begin{array}{cccc}
    \iota ': & \C ^{*} & \to  & \Pic ^{0}(C^{N})\\
    & t & \mapsto  & \O (p-\iota (t))
  \end{array}
  $$
  It follows that $\O(kp)\cong \O(k\iota(t))$ if and only if
  $\O(p-\iota (t))^{\otimes k}\cong \O $, i.e.~if $\iota'(t)$ is a
  $k$th root of unity.
\end{proof}

\subsection{The variety of \protect$L\protect$-osculating points}

In the setup of this section, the $L$-osculating points on the nodal
fibers can be used to define a global multi-section $\tilde{\sigma}
\subset \tilde{X}$. A detailed description of $\tilde{\sigma}$ will be
the key in our argumentation. For this, it is important to note that
the relative Picard group is locally divisible.

\begin{prop}\label{prop:local_divisibility}
  Let $k$ be the relative degree of the line bundle $L$, i.e.~the
  intersection number of $L$ with a fiber of $\pi$. Pick a point $\mu
  \in C$, fix a unit disk $\Delta \subset C$ centered about $\mu$ and
  set $X_{\Delta} := \pi^{-1}(\Delta)$. Then, after shrinking
  $\Delta$, if necessary, there exists a line bundle $L'\in \Pic
  (X_{\Delta})$ such that $kL'\cong L|_{X_{\Delta}}$.
\end{prop}

\begin{proof}
  As a first step, we will prove that $H^{2}(X_{\Delta},\Z )\cong \Z$.
  In order to see this, recall from deformation theory that ---after
  shrinking $\Delta$, if necessary--- $X_{\Delta}$ is of the form
  $$
  X_{\Delta }\cong \{(x,[y_0:y_1:y_2])\in \Delta \times \P_2\, |\,
  y_2y^2_1-y^3_0-f(x)y^2_0y_2\}
  $$
  where $f$ is a function $f\in \O(\Delta)$. In particular, if
  $N\subset X_{\Delta}$ is the non-normal locus and
  $\tilde{N}:=\eta^{-1}(N)$ its preimage in the normalization, then
  $N=\Delta \times \{[0:0:1]\}$ is a unit disc and $\tilde{N}$ either
  a unit disk or a union of irreducible components which are each
  isomorphic to $\Delta $ and meet in a single point. In this setup,
  we may use the Mayer-Vietoris sequence for reduced cohomology to
  calculate:
  $$
  \ldots \to \underbrace{H^{1}(\tilde{N},\Z )}_{=0}\to
  H^{2}(X_{\Delta },\Z )\to \underbrace{H^{2}(\tilde{X}_{\Delta },\Z
    )}_{=\Z }\oplus \underbrace{H^{2}(N,\Z )}_{=0}\to \ldots
  $$
  See \cite[prop.~3.A.7, p.~98]{BK82} for more information about
  the sequence.  Stefan~Helmke pointed out that $H^2(X_{\Delta},\Z)
  \cong \Z $ can also be shown by deforming $X_{\Delta}$ into a
  bundle of cuspidal plane cubics where the claim is obvious.
  
  Now choose a section $s\subset X_{\Delta }$ which is entirely
  supported on the smooth locus. After shrinking $\Delta $, this will
  always be possible.  Consider the exponential sequence
  \begin{equation} 
    \xymatrix{ {\ldots} \ar[r] & {H^1(X_\Delta, \O)}
      \ar[r]^{\alpha} & H^1(X_\Delta, \O^*) \ar[r]^{\beta} &
      {H^2(X_\Delta,\Z) \cong \Z } \ar[r] & {\ldots}} 
  \end{equation}
  The element $h:=(L|_{X_{\Delta }}-\O (ks))$ satisfies $\beta (h)=0$
  and is therefore contained in $\Pic^0(X_{\Delta})=\Image (\alpha)$.
  Let $h'\in \alpha ^{-1}(h)$ be a preimage and note, that, since
  $H^1(X_{\Delta},\O)$ is a $\C$-vector space, we can find an element
  $h''\in H^{1}(X_{\Delta },\O )$ such that $h'=k\cdot h''$ .  We may
  therefore finish by setting $L':=\alpha (h'')\otimes \O (s)$.
\end{proof}

The divisibility of $L$ implies that we can locally always find a
component of the osculating locus which is contained in the smooth
part $X_{\Reg} \subset X$.

\begin{cor}\label{cor:existence_of_special_section}
  Fix a point $\mu \in C$. If $\Delta \subset C$ is a sufficiently
  small unit disk about $\mu $, then there exists an $L$-osculating
  section $\sigma_1'\subset X_{\Delta}$ supported in the smooth locus
  of $X_{\Delta}$. More precisely, there exists a section $\sigma_1'
  \subset X_{\Delta ,\Reg }$ such that for all points $\mu \in \Delta
  $ the fiber $X_{\mu }:=\pi^{-1}(\mu )$ is either cuspidal or that
  the intersection $\sigma_1'\cap X_{\mu }$ is an $L$-osculating point
  of the fiber $X_{\mu}$.
\end{cor}

\begin{proof}
  Let $L'\in \Pic (X_{\Delta })$ be a line bundle such that $kL'\cong
  L$; the existence of $L'$ is guaranteed by
  proposition~\ref{prop:local_divisibility}.  Note that
  $R^{0}\pi_*(L')$ is locally free of rank one. Thus, after shrinking
  $\Delta $ if necessary, a section $s\in H^{0}(X_{\Delta },L')$
  exists whose restriction to any fiber of $\pi $ is not identically
  zero.  But since the relative degree of $L'$ is one, the restriction
  of $s$ to a fiber is a section which vanishes at exactly one smooth
  point of the fiber.  This point must therefore be $L$-osculating.
  Thus, the divisor $\sigma_1 \in \left| L'|_{X_{\Delta }}\right| $
  associated with the section $s$ contains only smooth $L$-osculating
  points and maps bijectively onto the base.
\end{proof}

This already gives a complete description of the $L$-osculating locus
in a neighborhood of a nodal fiber.

\begin{lem}\label{lem:spec_pts_at_nodal_pts}
  Let $C^{0}\subset C$ be the maximal (open) subset such that all
  fibers are nodal plane cubics. If $k$ is the degree of the
  restriction of $L$ to a fiber, then there exists a $k$-fold
  unbranched multisection $\sigma '\subset X_{\Reg }^0 :=
  \pi^{-1}(C^0)$ such that the restriction to any fiber $\sigma '\cap
  X_{\eta}$ is exactly the set of the $L$-osculating points of that
  fiber.
\end{lem}

\begin{proof}
  Let $\mu \in C^0$ be any point and $\Delta \subset C^0$ a small unit
  disk centered about $\mu $. The preimage $X_{\Delta} :=
  \pi^{-1}(\Delta) \subset X$ will then be isomorphic to $C^N \times
  \Delta $. Let $\sigma_1' \subset X_{\Delta,\Reg }$ be the
  $L$-osculating section whose existence is guaranteed by
  corollary~\ref{cor:existence_of_special_section} and find an
  isomorphism $\iota :\Delta \times \C ^{*}\to X_{\Delta ,\Reg }$ such
  that $\sigma_1' = \iota (\Delta \times \{1\})$. Apply
  lemma~\ref{lem:how_to_calculate_special_pts} to see that $\sigma '$
  is then given as
  $$
  \sigma =\{\iota (\Delta \times \{\xi _{i}\})\, |\, \xi
  ^{k}_{i}=1\}.
  $$
  Hence the claim.
\end{proof}

\begin{defn}
  Let $C^0\subset C$ be the maximal subset such that all fibers are
  nodal plane cubics. Let $\tilde{\sigma }\subset \tilde{X}$ be the
  closure of $\eta ^{-1}(\sigma ')$. We call the irreducible
  components $(\tilde{\sigma }_{i})_{i=1\ldots n}\subset
  \tilde{\sigma}$ the ``$L$-osculating (multi-)sections''.
\end{defn}

\subsection{\protect$L\protect$-osculating points in the neighborhood
  of a cusp}

As a next step we will find coordinates on $\tilde{X}_{\Delta} :=
\tilde{\pi}^{-1}(\Delta)$ where the $L$-osculating multisection
$\tilde{\sigma }$ can be written explicitly, even if $X_{\Delta }$
contains a cuspidal fiber.

\begin{prop}\label{prop:intersection}
  Assume that the preimage $\eta^{-1}(X_{\Sing})$ consists of two
  distinct sections $\tilde{\sigma}_0$ and $\tilde{\sigma}_{\infty}$
  and assume further that the $L$-osculating multisection
  $\tilde{\sigma} \subset \tilde{X}$ decomposes into $k$ irreducible
  components $(\tilde{\sigma}_i)_{i=1\ldots k}\subset X$ (which are
  then sections over $C$). If a point $p\in \tilde{\sigma }_0 \cap
  \tilde{\sigma}_{\infty}$ is given, then there is a unique index
  $1\leq j\leq k$ such that $p\not \in \tilde{\sigma}_j$ . All other
  components components $\tilde{\sigma}_a$, $\tilde{\sigma}_b$ with
  $a$, $b\not \in \{0,j,\infty \}$ do contain $p$ . If
  $m=\mult_p(\tilde{\sigma}_0, \tilde{\sigma}_{\infty})$ is the local
  intersection multiplicity of $\tilde{\sigma}_0$ and
  $\tilde{\sigma }_{\infty }$ at $p$, then
  $$
  \mult_p(\tilde{\sigma}_a,\tilde{\sigma}_b) =
  \mult_p(\tilde{\sigma}_{a},\tilde{\sigma}_0) =
  \mult_p(\tilde{\sigma}_a,\tilde{\sigma}_{\infty })=m.
  $$
\end{prop}

\begin{proof}
  Choose a unit disk $\Delta \subset C$ centered about $\mu :=
  \tilde{\pi }(p)$ and equip $\Delta $ with a coordinate $x$. Then
  $X_{\mu }:=\pi ^{-1}(\mu )$ is a cuspidal curve. After shrinking
  $\Delta $ we may assume that $\mu $ is the only point in $\Delta $
  whose preimage is cuspidal. By
  corollary~\ref{cor:existence_of_special_section}, we can find an
  index $j$ such that $p\not \in \tilde{\sigma}_j$. Therefore we can
  choose a bundle coordinate on $\tilde{X}_{\Delta }\cong \Delta
  \times \P_1$ so that we can write
  $$
  \begin{array}{rcl}
    \tilde{\sigma }_0 & = & \{([y_1:y_2],x)\in \P_1 \times \Delta \,
    |\, y_1 = x^my_2\}\\
    \tilde{\sigma }_{\infty } & = & \{([y_1:y_2],x)\in \P_1 \times
    \Delta \, |\, y_1 = -x^my_2\}\\
    \tilde{\sigma }_j & = & \{([y_1:y_2],x)\in \P_1\times \Delta \,
    |\, y_2 = 0\}
  \end{array}
  $$
  for an integer $m>0$. For a point $\nu \in \Delta $, $\nu \not
  =\mu$, the map $\eta \circ \iota_{\nu}$ with
  $$
  \begin{array}{cccc}
    \iota _{\nu }: & \C ^{*} & \to  & \tilde{\pi }^{-1}(x)\\
    & t & \mapsto  & [\nu ^{m}(t+1):(1-t)]
  \end{array}
  $$
  parameterizes the smooth part of $\pi ^{-1}(\mu )$. Apply
  lemma~\ref{lem:how_to_calculate_special_pts} with $\iota :=\eta
  \circ \iota_{\nu}$ and write
  $$
  \sigma =\{([y_{1}:y_{2}],x)\in \P _{1}\times \Delta \, |\, y_{1}(\xi
  -1)=x^{m}y_{2}(\xi +1),\, \xi ^{k}=1\}.
  $$
  The claim follows.
\end{proof}

The precise description of proposition~\ref{prop:intersection}
immediately implies that a projective bundle can be transformed into a
trivial bundle in a very simple manner.

\begin{cor}\label{cor:trafo_to_trivial}
  Under the assumptions of proposition~\ref{prop:intersection} above,
  if we define
  $$
  \tilde{D}_0 := \sum_{p\in \tilde{\sigma }_0 \cap
    \tilde{\sigma}_{\infty}} \mult_p(\tilde{\sigma}_0,
  \tilde{\sigma}_{\infty }) \cdot \tilde{\pi}(p)
  $$
  then $elt_{(\tilde{\sigma}_0,\tilde{D}_0)}: X\dasharrow \hat{X}$
  defines a map to the trivial bundle $\hat{X}\cong \P_1 \times C$.
  There exists a coordinate on $\P_1$ such that the strict transforms
  $\sigma_0$, $\sigma_{\infty}$ and $\sigma_i$ of $\tilde{\sigma}_0$,
  $\tilde{\sigma}_{\infty}$ and $\tilde{\sigma}_i$ are of the form
  $$
  \begin{array}{rcl}
    \sigma _{0} & = & \{[0,1]\}\times C\\
    \sigma _{\infty } & = & \{[1,0]\}\times C\\
    \sigma _{i} & = & \{[\xi _{i},1]\}\times C
  \end{array}
  $$
  where the $\xi_i$ are roots of unity.
\end{cor}

\begin{proof}
  It follows immediately from the construction of the elementary
  transformation that the intersection number between the strict
  transforms of $\tilde{\sigma}_0$ and $\tilde{\sigma}_{\infty}$ drops
  exactly by one with each transformation in the sequence
  $elt_{(\tilde{\sigma}_0,\tilde{D}_0)}$. In particular, the strict
  transforms $\sigma_0$, $\sigma_{\infty} \subset \hat{X}$ of
  $\tilde{\sigma}_0$ and $\tilde{\sigma}_{\infty}$ are disjoint.
  Likewise, it follows from proposition~\ref{prop:intersection} that
  the strict transforms $\sigma_i \subset \hat{X}$ of the
  $\tilde{\sigma}_i$ are sections which are mutually disjoint and
  disjoint from both $\sigma_0$ and $\sigma_{\infty}$. It follows that
  $\hat{X}$ must be the trivial bundle $\hat{X}\cong \P_1 \times C$
  and that $\sigma_0$, $\sigma_{\infty}$ and $\sigma_i$ are fibers of
  the projection $\pi_2:\hat{X}\to \P_1$.
  
  It remains to find the right coordinate on $\P_1$. To accomplish
  this, choose a general point $\mu \in C$. The fiber $X_{\mu} :=
  \pi^{-1}(\mu)$ will then be a nodal curve and there exists a point
  $p \in X_{\mu}$ such that $\O_{X_{\mu}}(kp)\cong L|_{X_{\mu}}$. It
  follows directly from lemma~\ref{lem:how_to_calculate_special_pts}
  that we find coordinates on $\tilde{X}_{\mu} = \eta^{-1}(X_{\mu})$
  such that $\tilde{\sigma}_0 \cap \tilde{X}_{\mu}$ corresponds to
  $[0,1]$, $\tilde{\sigma}_{\infty} \cap \tilde{X}_{\mu}$ corresponds
  to $[1,0]$, and the $L$-osculating points correspond to $[\xi_i:1]$
  where $\xi_i$ are roots of unity. Note that $elt_{(\tilde{\sigma}_0,
    \tilde{D}_0)}$ is isomorphic in a neighborhood of
  $\tilde{X}_{\mu}$ and use the coordinates on $\tilde{X}_{\mu}$ to
  obtain a global bundle coordinate on $\hat{X}\cong \P_1\times C
  \cong \tilde{X}_{\mu} \times C$. This coordinate will have the
  desired properties, and the proof is finished.
\end{proof}

\section{Proof of theorem~\ref{thm:main_characterization} }
\label{sec:Proof1}

Corollary~\ref{cor:trafo_to_trivial} already contains most arguments
needed for the proof of theorem~\ref{thm:main_characterization}.
Until the end of this section we use the notation of that theorem and
fix an ample bundle $L\in \Pic (X')$. Let $\tau :C\to C'$ be a base
change morphism which ensures that $\eta ^{-1}(X_{\Sing })$ consists
of two distinct section $\tilde{\sigma }_{0}$ and
$\tilde{\sigma}_{\infty}$, and that the $L$-osculating multisection
$\tilde{\sigma} \subset \tilde{X}$ decomposes into $k$ irreducible
components $(\tilde{\sigma}_i)_{i=1\ldots k}\subset X$ where $k$ is
the (positive) degree of $L$ if restricted to a $\pi $-fiber.
Consider the map $elt_{(\tilde{\sigma}_0,\tilde{D}_{0})}: X \dasharrow
\hat{X}\cong \P_1\times C$ which is defined in
corollary~\ref{cor:trafo_to_trivial} and note that we are finished if
we show that the inverse transformation
$elt_{(\tilde{\sigma}_0,\tilde{D}_{0})}^{-1}$ is of the form
$elt_{(\sigma_i,D_i)_{i=1\ldots n}}$ for effective divisors $D_i \in
\Div (C)$ with disjoint supports (here the $\sigma_i$ are defined as
in corollary~\ref{cor:trafo_to_trivial}). This, however, follows
directly from corollary~\ref{cor:inverse_trafo}. Actually, the
corollary shows that 
$$
D_i := \sum _{p\in \tilde{\sigma }_0 \cap
  \tilde{\sigma}_{\infty},\, \tilde{\sigma}_i\not \ni
  p}\mult_p(\tilde{\sigma }_{0}\cap \tilde{\sigma}_{\infty})\cdot
\tilde{\pi}(p),
$$
yield the desired transformation. This ends the proof of
theorem~\ref{thm:main_characterization}.

\section{Proof of theorem~\ref{thm:construction}}
\label{sec:Proof2}

\subsection{Sufficiency}

To begin the proof, we assume that $\sigma_i$ and $D_i$ are given as
in theorem~\ref{thm:construction} and that there exists a coordinate
on $\P_1$ such that the $\sigma_i$ correspond to roots of unity. We
will then show that the variety $X$ is projective. More precisely, we
fix an index $1\leq i\leq n$ and let $\tilde{\sigma}_i \subset
\tilde{X}$ be the strict transform of $\sigma _{i}$. We will show that
the image $\sigma_i' :=
\gamma_{(\tilde{\sigma}_0,\tilde{\sigma}_{\infty})}(\tilde{\sigma}_i)
\subset X$ is a $\Q$-Cartier divisor. Thus, a suitable multiple of
$\sigma_i'$ generates a relatively ample line bundle, and we are done.

If the construction of theorem~\ref{thm:construction} involves only
three sections $\sigma_0$, $\sigma_1$ and $\sigma_{\infty}$, then it
is clear that the strict transform $\tilde{\sigma}_1 \subset
\tilde{X}$ of $\sigma_1$ is disjoint from the strict transforms
$\tilde{\sigma}_0$ and $\tilde{\sigma}_{\infty}$. Thus, the image
$\sigma_1'$ does not meet the singular locus $X_{\Sing}$ of $X$ and is
therefore Cartier.

If the construction uses more than three sections and $\sigma_i'$ is
not already Cartier, let $\mu \in C$ be a point such that $\sigma_i'$
meets the singular locus $X_{\Sing}$ over $\mu$. By construction,
there exists a unique index $j$ such that $\mu \in |D_j|$. It follows
that the strict transform $\tilde{\sigma}_j$ of $\sigma_j$ does not
intersect $\tilde{\sigma}_0$ or $\tilde{\sigma}_{\infty}$ over
$\mu $:
$$
\tilde{\pi }(\tilde{\sigma }_{0}\cap \tilde{\sigma }_{j})\not \ni \mu
\quad \textrm{and}\quad \tilde{\pi }(\tilde{\sigma }_{\infty }\cap
\tilde{\sigma }_{j})\not \ni \mu .
$$
We can therefore find a suitable unit disc $\Delta \subset C$
centered about $p$ and we can find coordinates $x$ on $\Delta$ and a
bundle coordinate $[y_0:y_1]$ such that we can write
$$
\begin{array}{rcl}
  \tilde{\sigma }_{0} & = & \{([y_{0}:y_{1}],x)\in \P _{1}\times
  \Delta \, |\, y_{0}=x^{m}y_{1}\}\\
  \tilde{\sigma }_{\infty } & = & \{([y_{0}:y_{1}],x)\in \P _{1}\times
  \Delta \, |\, y_{0}=-x^{m}y_{1}\}\\
  \tilde{\sigma }_{j} & = & \{([y_{0}:y_{1}],x)\in \P _{1}\times
  \Delta \, |\, y_{1}=0\}
\end{array}
$$
An elementary calculation, using
lemma~\ref{lem:how_to_calculate_special_pts} and the assumption that
there exist coordinates where $\sigma_i$ and $\sigma_j$ are of the
form $\{\textrm{Root of unity}\}\times C$ shows that
$$
\begin{array}{rcl}
  \tilde{\sigma }_{i} & = & \{([y_{0}:y_{1}],x)\in \P _{1}\times
  \Delta \, |\, y_{0}=-\frac{\xi +1}{\xi -1}x^{m}y_{1}\}\\ 
  & = & \{([y_{0}:y_{1}],x)\in \P _{1}\times \Delta \, |\,
  \underbrace{y_{0}(\xi -1)+x^{m}y_{1}(\xi
    +1)}_{=:f(x,y_{0},y_{1})}=0\} 
\end{array}
$$
where $\xi$ is a root of unity. We fix a number $k$ such that
$\xi^k = 1$ and we will show that $\sigma_i'|_{X_{\Delta}}$ is a
$k$-Cartier divisor, i.e.~$k\cdot\sigma_i'|_{X_{\Delta}}$ is Cartier.
Recall from section~\ref{sec:Reduct_to_ruled} that the map $\gamma$
is locally given as
$$
\begin{array}{cccc}
  \gamma _{\Delta }: & \Delta \times \C  & \to  & \Delta \times \C
  ^{2}\\ 
  & (x,y_{0}) & \mapsto  & \left( x,y_{0}^{2}-x^{2m},\,
    y_{0}(y_{0}^{2}-x^{2m})\right) 
\end{array}
$$
In particular, the image of $\gamma_{\Delta }$ is isomorphic to
$\Spec R$, where $\gamma_{\Delta}^{\#}(R)\subset k[x,y_{0}]$ is the
subring generated by the constants $\C $, by $x$, by $y^2_0$ and by
the ideal $(y_0^2-x^{2m})$; see \cite[defn.~on p.~72]{Ha77} for the
notion of $\gamma_{\Delta}^{\#}$. Thus, to show that
$\gamma_{(\tilde{\sigma}_0,\tilde{\sigma}_{\infty})}(\tilde{\sigma}_1)$
is $k$-Cartier, it suffices to show that $f(x,y_{0},1)^{k}\in
\gamma_{\Delta}^{\#}(R)$.  We decompose $f^{k}$ as follows.
$$
\begin{array}{rcl}
  f(x,y_{0},1)^{k} & = & \sum _{i=0\ldots k}{k\choose
    i}(x^{m}-y_{0})^{i}(x^{m}+y_{0})^{k-i}\xi ^{k-i}\\
  & = &
  (x^{m}+y_{0})^{k}+(x^{m}-y_{0})^{k}+(x^{m}-y_{0})(x^{m}+y_{0})
  (\textrm{rest})\\
  & = & \sum _{i=0\ldots k}{k\choose
    i}[x^{m(k-i)}y_{0}^{i}+x^{m(k-i)}(-y_{0})^{i}] +
  (x^{m}-y_{0})(x^{m}+y_{0})(\textrm{rest})\\ 
  & = & \underbrace{2\cdot \sum _{i=0\ldots k,\, i\,
      \textrm{even}}{k\choose
      i}x^{m(k-i)}y^{i}_{0}}_{=:A}-
  \underbrace{(y^{2}_{0}-x^{2m})(\textrm{rest})}_{=:B}
\end{array}
$$
It is clear each summand of $A$ is in $\gamma_{\Delta }^{\#}(R)$
because it involves only even powers of $y_0$. Likewise, $B\in
\gamma_{\Delta}^{\#}(R)$ as $B$ is contained in the ideal
$(y^2_0-x^{2m})$. It follows that $f^k \in \gamma_{\Delta}^{\#}(R)$,
and we are done.

\subsection{Necessity}

It remains to show that the conditions spelled out in
theorem~\ref{thm:construction} are also necessary. For this assume
that $X$ is projective. We are finished if we can show that this
implies the existence of a coordinate on $\P_1$ such that $\sigma_0$,
$\sigma_{\infty }$ and $\sigma_i$ correspond to $[0,1]$, $[1,0]$ and
$[\xi_i,1]$ for certain roots of unity $\xi_i$. By
corollary~\ref{cor:trafo_to_trivial}, such a coordinate can be found
if the birational map $elt_{(\tilde{\sigma}_0,\tilde{D}_0)}: X
\dasharrow \hat{X}$ which was defined in
corollary~\ref{cor:trafo_to_trivial} is the inverse of the map
$elt_{(\sigma_i,D_i)} : X \dasharrow \tilde{X}$ which was used in
theorem~\ref{thm:construction} in order to construct the bundle $X$.
This, however, is exactly the statement of
lemma~\ref{lem:inverse_trafo}.


\begin{thebibliography}{Keb00}

\bibitem[BK82]{BK82}
G.~Barthel and L.~Kaup.
\newblock {\em Sur la Topologie des Surfaces complexes compactes}, chapter in
  Topologie des surfaces complexes compactes singuli\`eres, pages 61--297.
\newblock Number~80 in Semin. Math. Super. Les Presses de l'universit\'e de
  Montr\'eal, 1982.

\bibitem[Har77]{Ha77}
R.~Hartshorne.
\newblock {\em Algebraic Geometry}, volume~52 of {\em Graduate Texts in
  Mathematics}.
\newblock Springer, 1977.

\bibitem[Keb00]{Keb00a}
S.~Kebekus.
\newblock Families of singular rational curves.
\newblock LANL-Preprint math.AG/0004023, 2000.

\bibitem[Kol96]{K96}
J.~Koll\'{a}r.
\newblock {\em Rational Curves on Algebraic Varieties}, volume~32 of {\em
  Ergebnisse der Mathematik und ihrer Grenzgebiete 3. Folge}.
\newblock Springer, 1996.

\end{thebibliography}
\end{document}